# Heuristic algorithms for the operator-based relocation problem in one-way electric carsharing systems

Maurizio Bruglieri[a,*], Ferdinando Pezzella[b], Ornella Pisacane[c]

[a] Dipartimento di Design, Politecnico di Milano, Milano, Italy
[b] Dipartimento di Ingegneria dell'Informazione, Università Politecnica delle Marche, Ancona, Italy
[c] Facoltà di Ingegneria, Università degli Studi e-Campus, Novedrate (Como), Italy

**Abstract**

This paper addresses an Electric Vehicle Relocation Problem (E-VReP), in one-way carsharing systems, based on operators who move through folding bicycles between a delivery request and one of pickup. In order to deal with its economical sustainability, a revenue associated with each relocation request satisfied and a cost due to each operator used are introduced. The new optimization objective maximizes the total profit. To overcome the drawback due to the high CPU time required by the Mixed Integer Linear Programming formulation of the E-VReP, four heuristics, also based on general properties of the feasible solutions, are designed. Their effectiveness is tested on two sets of realistic instances. In the first one, all the requests have the same revenue. In the second one, the revenue of each request has a variable component related to the user's rent-time and a fixed one related to the customer satisfaction. Finally, a sensitivity analysis is carried out on both the number of requests and the fixed revenue component.

## 1 Introduction and literature review

Carsharing consists in the shared use of cars made available under the payment of a fare according to the time of use. It differs from car rental service since it provides for the use of the car only for a short time (usually for a fraction of an hour) in order to favor the sharing of each car among different people during the same day. Due to environmental sustainability issue, nowadays, several carsharing companies are providing the customers with the use of Electric Vehicles (EVs) rather than of traditional internal combustion engine cars. Real-life cases of electric carsharing systems are given by: Car2Go (http://www.car2go.com/) in Amsterdam and Vienna, for example; Autolib (www.autolib.eu) in Paris; Autobleue (http://www.auto-bleue.org) in Nice-Côte d'Azur (France).

* Corresponding author. Tel: +39 02.23995906. E-mail address: maurizio.bruglieri@polimi.it

The planning and the management of a carsharing service poses some important decision problems. About the planning, two significant issues concern the sizing of the fleet and the location of the parking stations (Cucu et al., 2009; Correia and Antunes, 2012; Jorge and Correia, 2013). About the management, some carsharing services permit one-way trips, which allow the user to pick up the vehicle in one station, and return it in another one. The flexibility offered by a one-way system makes it more attractive to users. But, on the contrary, it is harder to be managed since it involves a possible unbalancing between the demand and the availability of vehicles or vice versa between the request for returning the vehicles and the availability of vacant parking lots, making necessary a vehicle relocation. In such cases, the service provider has to develop strategies to relocate the vehicles and restore an optimal distribution of the carsharing fleet. Such strategies depend also on the available data and the main objective of the relocation. Barth and Todd (1999) propose the classification in: *static relocation*, based on the immediate needs of a particular parking lot; *historical predictive relocation*, based on an estimation of the requests made using historical data of the service or techniques of travel demand estimation; *exact predictive relocation*, based on the perfect knowledge of the request (like in the case of a carsharing service on reservation).

The vehicle relocation can be carried out by the user herself or by the service provider (Barth et al., 2004). In the first case, the user is motivated to car pool or to choose another parking station or reservation time (generally through the pricing lever); in the second case, the vehicles are moved by personnel either by trucks or directly driving each of them.

However, in general, operating with EVs rather than the traditional internal combustion engine ones complicates a lot the management of such systems as shown in Touati-Moungla and Jost (2012). In particular, the authors revise several problems related to the EV management (e.g., the poor battery range) from an optimization point of view. The carsharing service offered in the same location is studied also by Hafez et al. (2011), which, among the other things, minimize, the total travel time of relocation, through three different heuristics. Jung et al. (2014) address the problem of locating infrastructure for EVs (i.e., electric taxis in their case) by proposing a new model where the passenger demand is not known *a-priori* leading to a stochastic dynamic itinerary for each EV.

Kek et al. (2009) design a system based on a three-step optimization-trend-simulation for supporting carsharing operators in relocating the vehicles. Such a system is tested considering a realistic scenario of a carsharing company in Singapore.

Di Febbraro et al. (2012) model the complex dynamics of the carsharing system, using a discrete event system simulation. The paper considers the relocation made by both users and staff, and has a twofold objective: reducing the number of required staff and minimizing the number of carsharing vehicles to satisfy the system demand.

Correia et al. (2014) study the flexibility of the one-way carsharing systems and they propose a new mathematical formulation including the possibility that a user selects the station according to its vehicle availability and not only to the distance from her origin/destination. A real life case study, in Lisbon, is taken into account during the experimental campaign.

Very recent studies have been carried out in the work of Nourinejad and Roorda (2014) where the authors propose both a-priori benchmark model and a dynamic vehicle relocation optimization model. The latter is solved via a discrete-event simulator where the arrival of a user request constitutes an event. The numerical results discussed have shown that the optimization-simulation based model is suitable to determine a trade-off between the vehicle relocation time and the fleet size.

For the EV relocation problem, Bruglieri et al. 2014(a), propose therefore the use of a staff of workers that move easily and in eco-sustainable way from a delivery point to a pickup point by means of a folding bicycle that can be loaded in the trunk of the EV which needs to be moved. Such a new relocation approach generates a challenging pickup and delivery problem with features that have been never considered in the literature. We refer to such a problem as the **E**lectric **V**ehicle **Re**location **P**roblem (E-VReP )[1]. E-VReP shares some features with the 1-skip vehicle routing problem (Archetti and Speranza, 2004) and with the rollon-rolloff problem (Aringhieri et al., 2004). In particular, in all the three problems, only one item at time is assumed to be picked up and delivered. Moreover, the routes have to start from and end in a single depot without exceeding a given maximum time duration.

However, E-VReP is more challenging than the above mentioned problems since the distance covered by a vehicle depends also on the item picked up, i.e., the residual electrical charge of the EV picked up. In Bruglieri et al. (2014a), a Mixed Integer Linear Programming (MILP) formulation of the E-VReP is provided together with some techniques to speed up its solution.

We remark that in the E-VReP the set of delivery requests and of pickup requests are supposed to be known in advance and therefore, they are considered as input of the problem. This assumption is realistic because in the case of a carsharing system on reservation, it is possible to know *a priori* the real needs of EV relocation according to the fleet size, the parking size and the reserved carsharing demand. On the contrary (i.e. without reservations), the demand can be forecasted through different models and techniques proposed in the literature and mainly derived from studies in Logistics. For instance, Cucu et al. (2010) study a forecasting model to exploit customers' preferences so as to anticipate their needs and relocate the vehicles accordingly. Wang et al. (2010) propose an aggregated model at the station level for forecasting the total number of vehicles rented out and returned over time at each station.

In Bruglieri et al. 2014(b) real world alike instances of the E-VReP are created through a simulator applied to some data yielded by the Milano transport agency AMAT and by the main energy supplier company in Milano, A2A (www.a2a.eu ).

The aim of our paper is completely different from Bruglieri et al. 2014(a) since we want to investigate the economical sustainability of the EV relocation approach introduced in the latter where this aspect was neglected. This is important in order to understand the practicability of the E-VReP especially from the carsharing managers' point of view. To deal with the economical sustainability, we introduce the costs related to the use of the operators and a revenue associated with each relocation request satisfied. While, the original problem aims to handle as many requests as possible (neglecting the worker costs), the new problem, defined in Section 2 of this paper, wants to maximize the total profit.

---

[1] We have changed the acronym of the problem from E-VRP (given in the previous works) to E-VReP to avoid ambiguity with the Electric Vehicle Routing Problem.

However, this new problem reveals to be more challenging than the original E-VReP, and therefore, we design four heuristic approaches. Heuristic approaches were not in depth investigated in the previous works related to the E-VReP and then, they constitute an innovative contribution of this paper. The first two heuristics are greedy algorithms that consider the next request to be served according to the criterion of the nearest neighborhood and of the highest urgency, respectively (Section 3). The third algorithm is a more structured heuristic that iteratively builds the solution inserting a pair of compatible (pickup and delivery) requests at the time (Section 4). The policy followed to couple a pickup request with a delivery one consists in the minimum distance, while the priority order used to insert the next pair of requests in the route is guided by the so-called *Critical Factor* which is an indicator of the difficulty to serve a request if further delayed. At last the position where the pair of requests is inserted in the current route, is determined exploiting a proposition that guarantees when an insertion is feasible and another one that allows a-priori to compute the time extension due to such an insertion. The proofs of such propositions represent also original contributions of this paper. The forth heuristic is a randomized version of the third one. We test the effectiveness of the four heuristics on a set of benchmark instances, already proposed in literature, inspired to a real word context (in the city of Milano). These results are compared to those provided by two MILP formulations (Section 5). Moreover, to test the economical sustainability issue, we build a new set of instances where the revenue associated with each relocation request depends on two components: a variable one, proportional to the rent-time of the request (i.e., the time in which the customer associated with the request uses the carsharing service) and a fixed one. The latter allows modeling a "future revenue" related to the customer satisfaction since a satisfied customer will (likely) ask for the service also in the future.

To experiment with the economical sustainability, we also perform a sensitivity analysis varying the two main input parameters of the E-VReP, i.e., the number of requests (instance size) and the fixed revenue component associated with each relocation request (Section 6). Finally, we draw some conclusions and future works (Section 7).

## 2    E-VReP definition and extension

We recall the description and the notation already used in literature for the E-VReP. A one-way carsharing service with a homogeneous fleet of EVs is given. Let $L$ be the maximum distance that an EV can cover when its battery is fully charged. Such a distance depends on the kind of EV considered (in the experimental campaign we assume that $L$ =150 km). When the battery of an EV is not fully charged, the maximum distance that can be covered is assumed to be linearly proportional to the residual charge of the battery (i.e. an EV with residual charge at 50% can cover $L/2$ km). Concerning the battery recharging at a station, we can assume that the battery charge level increases linearly with the time. Let $\Gamma$ indicate the time to fully recharge a com-pletely exhausted battery. We suppose that every EV is always picked up and returned to a parking station lot equipped with a recharge dock so as it is recharged when it is not used. Since in a one-way carsharing service the cars can be returned to parking stations different from those where they are picked up, some of them need to be moved in order to prevent that a station runs out of either EVs or parking lots.

Let *D* be the set of delivery requests (i.e., the requests of EVs to deliver in order to prevent that a station runs out of EVs) and let *P* be the set of pickup requests (i.e. requests of EVs that need to be moved to vacant parking lots). Each relocation request $r \in PUD$ is characterized by a parking location $v_r$, i.e., a node of the road network, by the residual charge of the battery $\rho_r$ and by a time window $[\tau_r^{min}, \tau_r^{max}]$ where $\tau_r^{min}$ and $\tau_r^{max}$ represent respectively the earliest time and the latest time when is allowed to carry out the request $r$. For instance if $r$ is a pickup request then $\tau_r^{min}$ is the time before which the EV is not available while $\tau_r^{max}$ is the time after which picking up the EV is not convenient (since from $\tau_r^{max}$ it may be used by some other customer). For a delivery request $r$, $\rho_r$ indicates the minimum charge level that the EV battery must have at time $\tau_r^{max}$. Therefore, if an EV is delivered before $\tau_r^{max}$, the charge level of its battery may be less than $\rho_r$. This is allowed on condition that at least $\rho_r$ is achieved at $\tau_r^{max}$, considering that its battery is recharged after the delivery. Whereas for a pickup request $r$, $\rho_r$ indicates the battery charge level at $\tau_r^{min}$. Since the fleet of EV is homogeneous each delivery request can be satisfied picking up every EV of a pickup request on condition that it is compatible for time windows and battery charge level. Given a team of *K* workers which leave a single depot, even at different times, using folding bicycles, we want to determine their routes and schedules in such a way that the objective to maximize is the number of satisfied requests. Each route consists of an alternating sequence of pickup requests and delivery requests, such that its duration does not exceed a given threshold *T* (i.e., the duty time of the workers), it ends at the depot and both the time windows and battery charge level constraints are respected. This is the original statement of the E-VReP.

To deal with the economical sustainability of the E-VReP, we change its objective assuming that a revenue $rev_i$ is associated with each relocation request *i* satisfied and a cost *C* is associated with each worker used. Thus, the E-VReP objective is modified into the maximization of the total profit given by the difference between the total revenue, represented by the sum of the revenues of all the relocation requests satisfied, and the total cost, obtained multiplying by *C* the number of the workers used in the fixed time horizon. The recharging cost is not directly considered because it can be assumed fixed. In fact, the carsharing companies usually pay a fixed monthly fee to the electric energy providers that is proportional to their fleet size (being this a measure of the electricity consumption). Moreover, through the new objective, we obtain an extension of the original E-VReP since the previous objective can be obtained again setting $C = 0$ and $rev_i = 1$ for every relocation request $i \in PUD$. All the input parameters used by the E-VReP are summarized in Table 1.

| Symbol | Meaning |
|---|---|
| $P$ | pickup set |
| $D$ | delivery set |
| $R=P\cup D$ | set of requests (either pickup or delivery type) |
| $v_r$ | parking location of request $r$ |
| $[\tau_r^{min}, \tau_r^{max}]$ | time window of request $r$ |
| $\rho_r$ | battery level demanded by request $r$ |
| $rev_r$ | revenue of request $r$ |
| $l_{ij}$ | length of the shortest path from $v_i$ to $v_j$ |
| $s''$ | average speed of a bike |
| $s'$ | average speed of an EV |
| $q'$ | average time for parking an EV and taking the bike up from the EV trunk |
| $q''$ | average time for unloading the bike from the EV trunk and leaving the parking lot with the EV |
| $\Gamma$ | time for a full battery recharge |
| $L$ | maximum distance travelled by a fully charged EV |
| $0$ | depot |
| $K$ | number of available workers |
| $T$ | duty time of each worker |
| $C$ | worker cost |

**Table 1** Input parameters of the E-VReP

An instance of the problem with five delivery requests and five of pickup is depicted in Fig. 1. In this example, it is assumed $C = 30$ €, a duty time equal to 4 hours for each worker, a revenue of 10 € for each request satisfied, a travel time (by EV) between each pair pickup-delivery of 20 minutes and a travel time (by bike) between each pair delivery-pickup of 30 minutes. Fig. 1 shows the optimal solution of the original E-VReP where the number of workers used is 2 and therefore, the total cost is equal to 60 €. Since six requests are satisfied, the total profit is zero. This behavior is due to the fact that the original E-VReP aims to maximize the number of requests served. While, the optimal solution of the new problem introduced in this paper consists of only the route performed by the operator $W_1$. Although fewer requests are satisfied, the total profit is not zero but equal to 10 €.

Hereafter the MILP formulation of the E-VReP proposed by Bruglieri et al. 2014(a) is denoted as MILP1. To model the economical sustainability issue, we change the objective function of MILP1 in the following way:

$$\max \sum_{k=1}^{K} \sum_{(i,j) \in A: i \neq 0} rev_i x_{ijk} - C \sum_{k=1}^{K} \sum_{(0,j) \in A} x_{ojk} \qquad (1)$$

where $x_{ijk}$ represent the routing variables of the $k$-th worker, $A$ is the set of the arcs of the modeling graph, linking every pair of compatible relocation requests, and 0 represents the depot node. Hereafter, we denote this new formulation as MILP2.

For this extension of the E-VReP, in the next two sections, we propose innovative heuristics suitable to solve instances with a significant number of requests.

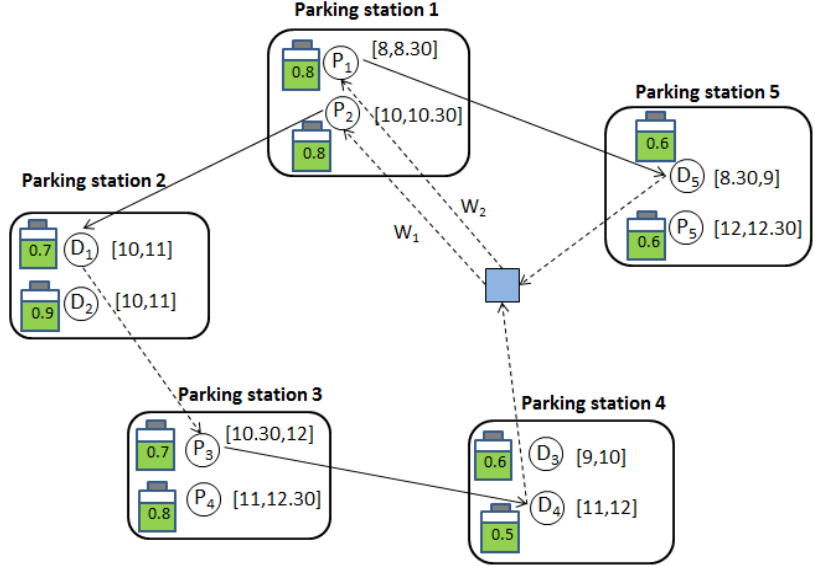

**Fig.1.** Instance of the E-VReP with five pickup requests ($P_i$ for $i = 1,...,5$) and five delivery requests ($D_i$ for $i = 1,...,5$). Both the battery charge level and time window $[\tau_r^{min}, \tau_r^{max}]$ are indicated beside each request. Two workers $W_1$, $W_2$ leave a common depot, indicated by the square node, to relocate the EVs. The dashed arcs denote that a worker is biking while the solid arcs denote that is driving an EV.

## 3 Two Greedy Heuristics for the E-VReP

In this section, two constructive heuristics are described. Both of them are outlined in Fig. 1. They differ for the way the next request to be served is selected (line 7 of Fig.1). In particular, in the *Nearest Neighborhood Heuristic* (NNH), the next request to be handled by a worker is the feasible one at minimum distance from the currently served request. Instead, in the *Most Urgent Heuristic* (MUH), the next request to be handled is the feasible one with the lowest right side time window i.e. the most urgent feasible one. Before introducing the feasibility conditions (indispensable to select the next request), the following clarifications are needed. In fact, in both the heuristics, the arrival time at a pickup request $p$, indicated with $a_p$, is settled as:

$$a_p = a_{d-1} + w_{d-1} + q' + \frac{l_{d-1p}}{s''} \qquad (2)$$

where *d-1* denotes the request of delivery previously served (for the remaining notation refer to Table 1). While, the arrival time at a delivery request *d* is fixed as

$$a_d = a_p + w_p + q'' + \frac{l_{pd}}{s'} \qquad (3)$$

where *p* is the request of pickup previously served.

When the request *p* is the first inserted, *d* represents the depot and $a_0$ is set as

$$a_0 = \tau_p^{min} - \frac{l_{0p}}{s''} \qquad (4)$$

The waiting time $w_d$ to satisfy a delivery request *d* is given by:

$$w_d = \max\{0, \tau_d^{min} - a_d - q'\} \qquad (5)$$

while the waiting time $w_p$ to satisfy a pickup request *p*, is given by:

$$w_p = \max\{0, \tau_p^{min} - a_p\} \qquad (6)$$

The slight asymmetry in the two kinds of waiting times is due to the fact that the delivery can start at $\tau_d^{min} - q'$ to deliver the EV exactly at $\tau_d^{min}$, while the EV cannot be picked up at $\tau_p^{min} - q''$ since it is not available yet.

In both NNH and MUH, the generic request *r* to handle is selected among the ones that satisfy the following <u>necessary and sufficient</u> feasibility conditions.

**Case 1**- The request previously handled by the worker is of delivery, then $r = p \in P$. The following two conditions have to hold:

$$a_p \leq \tau_p^{max} \qquad (7)$$

$$\max\{a_p, \tau_p^{min}\} + q'' + \min_{d \in \tilde{D}}\left\{\frac{l_{pd}}{s''} + \frac{l_{d0}}{s''}\right\} + q' - start^\omega \leq T \qquad (8)$$

where $\tilde{D}$ denotes the set of delivery requests compatible with *p* not yet served. The condition (7) ensures that the request *p* is satisfied within the latest allowed time and condition (8) guarantees that the duration of the route does not exceed the threshold *T*.

**Case 2**-The request previously handled by the worker is of pickup (let us indicate it by *p*), then $r = d \in D$. The following four conditions have to hold:

$$a_d \leq \tau_d^{max} \qquad (9)$$

$$\max\{a_d, \tau_d^{min}\} + q' + \frac{l_{d0}}{s''} - start^\omega \leq T \qquad (10)$$

$$\min\left\{\rho_p + \frac{a_p + w_p - \tau_p^{min}}{\Gamma}, 1\right\} - \frac{l_{pd}}{L} \geq 0 \qquad (11)$$

$$\min\left\{\rho_p + \frac{a_p + w_p - \tau_p^{min}}{\Gamma}, 1\right\} - \frac{l_{pd}}{L} + \frac{\tau_d^{max} - a_d}{\Gamma} \geq \rho_d \qquad (12)$$

The conditions (9) and (10) work like conditions (7) and (8), respectively. Condition (11) ensures that the current battery level is sufficient to go from $v_p$ to $v_d$ along the minimum path. Condition (12) guarantees that an EV is delivered with a battery level such that at the time $\tau_d^{max}$ a charge level not lower than $\rho_d$ is reached. If the size of the pickup set *P* is not equal to the size of the delivery set *D*, then both the greedy heuristics end as soon as no further element in the smallest of these two sets, can be inserted into a route. Thus, some requests belonging to the biggest of the two sets are sacrificed in greedy way.

**Greedy Heuristic (*K*, *R*, *T*, *Ω*):**
```
1.   Ω ≔ ∅
2.   While (|Ω| ≤ K and R ≠ ∅) do
3.      Open a new route ω;
4.      Ω ≔ Ω∪{ω};
5.      R' ≔ R;
6.      While (R' ≠ ∅) do
7.         Select a feasible request r ∈ R';
8.         R' ≔ R'\{r};
9.         If (time duration of ω∪{r} ≤ T)
10.           Add r to the end of ω;
11.           Set a_r and w_r;
12.           R ≔ R\{r};
13.        EndIf
14.     EndWhile
15.  EndWhile
16.  Return Ω.
```

**Fig. 2.** The general outline of the greedy heuristics

## 4   Two structured heuristics for the E-VReP

In this section we introduce two heuristics for the E-VReP smarter than those presented in Section 3 since they exploit some properties of the feasible solution set. The two heuristics are named *Critical Heuristic* (CH) and *Randomized Heuristic* (RH).
First of all, some requests have to be possibly removed to balance the test instance, i.e., in order to obtain |P|=|D|, since the relocation method here considered needs always to alternate a pickup request with a delivery one. Therefore, we start both the heuristics with a *preprocessing phase*. Then, in the *construction phase* the solution is iteratively built inserting a pair of compatible (pickup and delivery) requests at the time. For this purpose, it is necessary to define:
  i)     the policy followed to couple a pickup request with a delivery one;
  ii)    the priority order used to insert the next pair of requests in the route;
  iii)   the position in the current route where the pair of requests is inserted.
We now analyze more in depth all these elements.

First of all, we introduce to the so-called *Critical Factor* (*CF*) which is an indicator of the difficulty to serve a request based on its level of urgency. If the request is not coupled with those still available it may not be served any more. The *CF* is differently computed according to the type of request. For every pickup request *p*, *CF(p)* is defined as

$$CF(p) = \max_{d \in \widetilde{D}}\{\tau_d^{max} - t_{pd}\} - \tau_p^{min} \tag{13}$$

where $\widetilde{D}$ indicates the set of delivery requests compatible with $p$ not yet served and $t_{pd}$, the time necessary to go from $v_p$ to $v_d$. Since $\max_{d \in \widetilde{D}}\{\tau_d^{max} - t_{pd}\}$ represents the latest departure time from the location of request $p$ in order to serve a feasible (unserved) delivery request, the smaller the gap with $\tau_p^{min}$ is, the more urgently the request needs to be served. In fact, if the delivery request which the maximum is achieved with, is no longer available then the request may remain without the possibility to be coupled with a delivery one (since it becomes more unlikely to satisfy the time window constraint). Similarly, for a delivery request $d$, $CF(d)$ is defined as

$$CF(d) = \tau_d^{max} - \min_{p \in \widetilde{P}}\{\tau_p^{min} + t_{pd}\} \tag{14}$$

where $\widetilde{P}$ indicates the set of pickup requests compatible with $d$ not yet served.

The feasibility of requests $p$ and $d$ is checked through the following three <u>necessary</u> conditions:

$$\tau_p^{min} + t_{pd} + q'' + q' \leq \tau_d^{max} \tag{15}$$

$$\rho_p - \frac{l_{pd}}{L} + \frac{\tau_d^{max} - \tau_d^{min}}{\Gamma} \geq \rho_d \tag{16}$$

$$t_{0p} + t_{d0} + \max\{t_{pd} + q'', \tau_d^{min} - \tau_p^{max}\} + q' \leq T \tag{17}$$

The condition (15) is necessary to the respect of the time window of $d$ considering the departure from $v_p$ at $\tau_p^{min}$. Condition (16) is necessary to the battery level feasibility. In fact, if it is violated, then the battery level $\rho_d$ cannot be achieved at $\tau_p^{max}$ since the maximum time that can elapse between the pickup request and the delivery one is equal to $\tau_p^{max} - \tau_p^{min}$. Condition (17) is necessary in order to guarantee that the duty time $T$ of the operators is never exceeded. In this phase of the algorithm, we cannot guarantee that the pair of requests $p$ and $d$ is feasible because the arrival time and the waiting time at $v_p$ and $v_d$ have not been fixed yet. For this reason, we cannot check if they also satisfy the sufficient feasibility conditions (9)-(12).
The requests with a negative *CF* have to be removed since infeasible with regard to their time windows.
The *preprocessing phase* consists in the following steps for both the two heuristics:
1. the requests are sorted by increasing values of *CF*;
2. if |P|>|D| then the first |P|-|D| pickup requests are removed;
3. otherwise, if |D|>|P|, the first |D|-|P| delivery requests are removed.

The removed requests will be not handled since they are considered rejected.

Concerning items i) and ii) of the *construction phase*, at each iteration, the CH heuristic selects the most critical request ($p$ or $d$), i.e. the one with minimum value of *CF*, and couples it with the request ($d$ or $p$) satisfying the necessary feasible conditions (15)-(17), whose parking is at minimum distance. If any request cannot be coupled with the currently critical one, the latter will be removed.

Once the first pair $(p_1, d_1)$ to be inserted has been detected, the route $\omega$ is initialized with $\omega = \{0, p_1, d_1, 0\}$. Besides, the arrival and the waiting times of the requests $p_1$ and $d_1$ are initialized in the following way:

$$a_{p1} = min\{\tau_{p1}^{max}, \mu - t_{p1d1} - q' - q''\} \tag{18}$$

where

$$\mu = \max\{\tau_{d1}^{min}, \tau_{p1}^{min} + t_{p1d1} + q' + q''\} \tag{19}$$

and

$$w_{p1} = 0 \tag{20}$$
$$a_{d1} = a_{p1} + t_{p1d1} + q' + q'' \tag{21}$$

$$w_{d1} = \mu - a_{d1} \tag{22}$$

This particular setting guarantees time window feasibility and null waiting time at $v_{d1}$ if $\tau_{p1}^{max}$ is large enough i.e. if $\tau_{d1}^{min} - t_{p1d1} - q' - q'' \leq \tau_{p1}^{max}$. More in general, given a feasible route $\omega = \{0, p_1, d_1, p_2, d_2, \ldots, p_n, d_n, 0\}$ and a pair $(p,d)$ that has to be inserted between $d_{i-1}$ and $p_i$, the arrival time and the waiting time of $p$ and $d$ are set by the CH heuristic as in (2), (3), (5), (6).

The main idea beside these equations is to fix the arrival times to both $p_1$ and $d_1$ in order not to wait at each of them. With the aim to better clarify the meaning of the equations (18)—(22) we give the following numerical example.

**Example 1:** It is assumed that the pair $(p_1, d_1)$ is selected as the first to be inserted. The related time windows are [100, 200] and [180, 300], respectively, while, the time to reach $d_1$ from $p_1$ (by EV) is $t_{p_1 d_1} = 50$ and the travel time $t_{0p_1}$ (by bike) is 20. According to equation (20), the waiting time at $p_1$ is zero since the vehicle is assumed leaving the depot in order to reach it not before the left-side extreme of the time window (i.e., 100). In this way, the arrival time to $p_1$, according to equation (18), is $a_{p1} = min\{200, \mu - 50 - q' - q''\}$ where, according to (19), $\mu = \max\{180, 100 + 50 + q' + q''\} = 180$, assuming the two input parameters $q'$ and $q''$ are equal to 1. Then the arrival time to $p_1$ is $a_{p1} = min\{200, 128\} = 128$ and the arrival time to $d_1$, according to (21), is $a_{d1} = 128 + 50 + 2 = 180$. The latter is exactly equal to the left-side extreme of its time window and therefore, according to (22), the waiting time at $d_1$ is $w_{d1} = 180 - 180 = 0$. In conclusion, computing these times as in equations (18)—(22) guarantees not to wait at both the nodes $p_1$ and $d_1$. It is easy to see that if, instead, the driver leaves the depot for instance at 0, the waiting time at $d_1$ will be 26 rather than 0.

After the insertion of the first pair of requests, the next pair $(p_2, d_2)$ will be inserted either before or after $(p_1, d_1)$. Generally, at this step of the procedure (item *iii* of the *construction phase*), different positions for the insertion in the current route of the new pair of requests could be feasible. Among them we want to choose the one with the minimum time extension of the route. For this purpose we exploit the following two propositions that allow determining a-priori the time extension of the route due to the insertion of a couple of requests and establishing if the insertion can be done preserving the route feasibility, respectively. In this way, the computational times of the solution approach can be reduced.

**Proposition 1**

Given a feasible route $\omega = \{0, p_1, d_1, p_2, d_2, \ldots, p_n, d_n, 0\}$ and a pair $(p,d)$ that can be feasibly inserted between $d_{i-1}$ and $p_i$, the Time Extension (*TE*) of $\omega$ due to such an insertion is given by:

$$TE(\omega) = max\{0, max\{max\{a_p, \tau_p^{min}\} + q'' + t_{pd}, \tau_d^{min}\} + q' + t_{dp_i} - a_{p_i} - \sum_{k=p_i}^{d_n} w_k\} \tag{23}$$

**Proof**:

After the insertion of the pair $(p,d)$, the new arrival time to $p_i$ is equal to

$$max\{max\{a_p, \tau_p^{min}\} + q'' + t_{pd}, \tau_d^{min}\} + q' + t_{dp_i}.$$

Then, the time extension of the route is given by the difference between the new arrival time to $p_i$ and the one before the insertion (i.e., $a_{pi}$) reduced by the sum of the waiting times at all the nodes next to it. This exactly corresponds to the formula (23).
∎

**Proposition 2**

Given a feasible route for the E-VReP, $\omega = \{0, p_1, d_1, p_2, d_2, \ldots, p_n, d_n, 0\}$, <u>sufficient</u> conditions to insert a feasible pair of requests $(p,d)$ between $d_{i-1}$ and $p_i$ maintaining the feasibility are:

$$a_{d_{i-1}} + w_{d_{i-1}} + t_{d_{i-1}p} \leq \tau_p^{max} \tag{24}$$

$$max\{a_{d_{i-1}} + w_{d_{i-1}} + t_{d_{i-1}p}, \tau_p^{min}\} + t_{pd} + q' + q'' \leq \tau_d^{max} \tag{25}$$

$$min\{\rho_p + (max\{a_{d_{i-1}} + w_{d_{i-1}} + t_{d_{i-1}p}, \tau_p^{min}\} - \tau_p^{min})/\Gamma, 1\} - l_{pd}/L \geq 0 \tag{26}$$

$$min\{\rho_p + (max\{a_{d_{i-1}} + w_{d_{i-1}} + t_{d_{i-1}p}, \tau_p^{min}\} - \tau_p^{min})/\Gamma, 1\} - l_{pd}/L +$$
$$(\tau_d^{max} - max\{max\{a_{d_{i-1}} + w_{d_{i-1}} + t_{d_{i-1}p}, \tau_p^{min}\} + q'' + t_{pd}, \tau_d^{min}\})/\Gamma \geq \rho_d \tag{27}$$

$$max\{max\{a_{d_{i-1}} + w_{d_{i-1}} + t_{d_{i-1},p}, \tau_p^{min}\} + q'' + t_{pd}, \tau_d^{min}\} + q' + t_{dp_i} - a_{p_i} \leq min_{r=p_i,d_i,p_{i+1},\ldots,d_n}\{\sigma_r + \sum_{k=p_i}^{r} w_k\} \tag{28}$$

$$a_{d_n} + w_{d_n} + t_{d_n 0} - (a_{p_1} - t_{0p_1}) + TE(\omega) \leq T \tag{29}$$

where in (28) $\sigma_r = \tau_r^{max} - a_r - w_r$ represents the maximum postponement related to the request $r$ and in (29) $TE(\omega)$ indicates the time extension given by (23).

**Proof**:

The conditions aforementioned are sufficient because if we leave unchanged the arrival times and the waiting times of every request preceding $p$ and we set those of $p$ and $d$ according to (2),(3),(5),(6), the necessary and sufficient conditions for feasibility (7)-(12) hold for ($p,d$) and for every successive request (besides, of course, every request preceding $p$). The condition (7) holds for $p$ thanks to (24), condition (8) thanks to (29). While, condition (9) holds for $d$ thanks to (25), condition (10) thanks to (26), condition (11) thanks to (27) and finally condition (12) thanks to (29). In particular condition (27) guarantees that the starting battery level ($\rho_p$) increased by the recharging level obtained in the time elapsed from $\tau_p^{min}$ until the instant of pickup ($\max\{a_{d_{i-1}} + w_{d_{i-1}} + t_{d_{i-1}p}, \tau_p^{min}\}$), decreased by the battery consume for going from $v_p$ to $v_d$ ($l_{pd}/L$) and increased by the battery level reached in the time elapsed from the instant of delivery ( $\max\{\max\{a_{d_{i-1}} + w_{d_{i-1}} + t_{d_{i-1}p}, \tau_p^{min}\} + q'' + t_{pd}, \tau_d^{min}\}$) until the maximum time allowed for the delivery ($\tau_d^{max}$), has to be greater or equal than $\rho_d$.

The conditions (7) and (9) hold also for every pickup request and delivery request successive to $p$, respectively, thanks to condition (28). In particular, its left side represents the difference between the time needed to reach $p_i$ and the one requested before the insertion of the pair ($p,d$). In other words, it represents the time extension of the route in reaching $p_i$ due to the insertion of the new pair ($p,d$). While the right side of (28) represents the maximum allowed time extension taking into account the maximum postponement of each request $r$ (denoted by $\sigma_r$), successive to $d$ and their waiting times can allow reducing the time extension of the route.

Regarding the battery level feasibility of all the delivery requests successive to $d$, the fact that, after the insertion of the pair ($p,d$), their arrival times to the corresponding stations may be only postponed, assures that both conditions (10) and (11) continue to hold since the term $\min\{\rho_p + \frac{a_p + w_p - \tau_p^{min}}{\Gamma}, 1\}$ can only increase.

Finally the duty time feasibility of the route is guaranteed by condition (29). ∎

According to Proposition 1, if several insertions are feasible for a pair of requests ($p,d$), then the insertion minimizing the time extension given by (23) is chosen. A route is closed when no pair of requests can be handled within the duty time $T$. Until workers are still available, a new route is open for handling the remaining requests according to the same operations applied to the first route.

Concerning the heuristic RH, a candidate request is chosen randomly and then, the joinable request is selected according to the necessary feasibility conditions (15)-(17). The right insertion position is then detected again according to both the conditions (24)-(29) and the minimization of the time extension (23), as in CH. Since in this way different feasible solutions can be generated, the constructive phase is repeated several times and the best solution is selected as the final one according to the objective function considered. Therefore, the results obtained with this heuristic can be different according to the objective function considered, i.e., either the total number of served requests or the total profit. This behavior is only related to the RH while the other heuristics are equal for the two objective functions.

## 5 Computational results

In this section, the experimental campaign carried out by running all the four heuristic solution approaches proposed in this paper is described.
All the heuristics have been coded in Java, in Eclipse (open-source integrated software development environment). The maximum number of iterations has been set to 10,000 runs, for the heuristic RH. The MILP formulations have been implemented in AMPL (Fourer et al., 2002) and solved through the state of the art solver CPLEX12.5 with a CPU time limit of 43,200 seconds, for comparison purposes. The experiments have been run on an Intel® Core™ i7-2630QM CPU @ 2.00 GHz and 8 GB RAM.
In order to assess the performances of the heuristic approaches, they have been tested on two sets of instances: the first one (AMAT, for short) is the benchmark set of instances for the E-VReP described in Bruglieri et al. 2014(b) and characterized by 22 requests, on average while the second one (V-AMAT, for short) has been generated specifically in this work and its instances are characterized by a variable profit associated with each request as detailed in the following. Both of them have been generated by using a simulator, where the relocation requests have been estimated taking into account the origin-destination traffic matrix yielded by AMAT, the Milano transport agency (AMAT, 2005). While, the location and the capacity of the docking stations considered are those installed by A2A, the main energy supplier company in Milano (A2A, 2013).
The data provided by the AMAT agency concerns the private car movements and are represented through the Origin-Destination (O-D) matrix from/to different zones of Milano, with trips having different aims (business, study, occasional, etc) and in different time-slots of the day: morning (from 7.00 a.m. to 10.00 a.m.), not-peak (from 10.00 a.m. to 4.00 p.m.) and evening (from 4.00 p.m. to 8.00 p.m.). In particular, we consider the data regarding the occasional trip since it represents a more common situation of carsharing using.
Fig. 3 shows the location of the A2A charging stations (i.e., 5 stations with 4 slots and 21 with 2 slots). The AMAT O-D zones have been intersected with a circular boundary of 500 meters around each charging station, representing the area easily reachable by walking from the station. The intersection allows estimating the potential number of movements that could be carried out with carsharing service rather than with private car. Moreover, such values have been multiplied by 0.5% in order to consider only a real-like percentage of the potential demand. Such a percentage is consistent with the current usage of the Milano carsharing service.
With the aim of estimating the requests of the relocation staff, Bruglieri et al. 2014(b) evaluated the unbalances due to the estimated carsharing travel demand through a simulator. Such a carsharing simulator, developed in Matlab, has a time-step feed with the data about the station capacities, the travel times between pairs of stations, and the travel demands. More specifically, the simulator randomly generates the vehicle movements between pairs of stations following a probability distribution accord-ing to the carsharing travel demand estimated from the AMAT data. At each simulated minute, the inventory of each station and both the position and the charge level of the EVs are updated according to their movements.

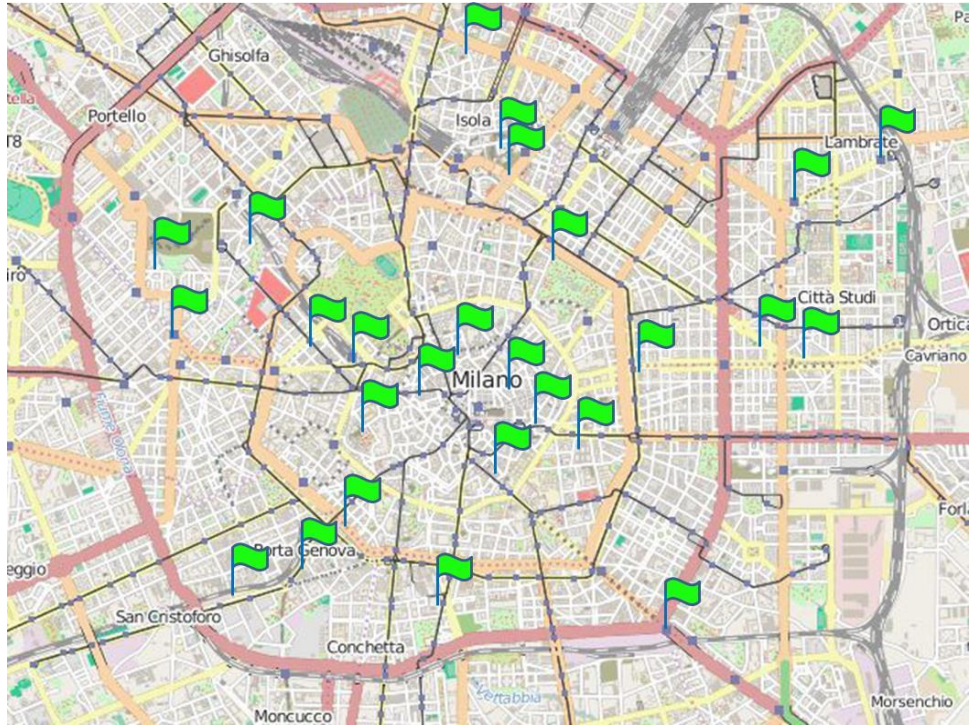

**Fig. 3.** Map of Milano with the locations of the A2A charging stations indicated by the flags (source: OpenStreetMap).

During the simulation, if a station *s* runs out of EVs and a user asks for an EV, then a delivery request *r* for the E-VReP is generated with: $v_r$ equal to the node of the road network corresponding to *s*, $\tau_r^{min}$ equal to the minute at which such an event happens and $\tau_r^{max}$ equal to the earliest minute when an EV arrives to *s* (or $\tau_r^{max}=\infty$ if the latter event never happens). Similarly, a pickup request for the E-VReP is generated as soon as a parking station becomes full and a user wants to give an EV back in such a station. Finally, $\rho_r$ is computed considering a starting random level of charge of the EV associated with the request *r* and decreasing or increasing such a value during the simulation according to the use of the EV or its parking time at the recharging station, respectively. In this way, the simulator returns the set of pickup requests *P*, the one of delivery requests *D* and for each request *r*, the values $\tau_r^{min}$, $\tau_r^{max}$ and $\rho_r$.
The other input data are summarized in Table 2.

| Input data | Value |
| --- | --- |
| $T$ | 300 minutes |
| $s'$ | 25 km/h |
| $s''$ | 15 km/h |
| $q'$ | 1 minute |
| $q''$ | 1 minute |
| $L$ | 150 km |
| $\Gamma$ | 240 minutes |
| $C$ | 60€ |

**Table 1.** Values of the main input data

Concerning the AMAT set of instances, the testing phase has been carried out considering two different scenarios. In the former, the worker cost is not taken into account during the optimization process and therefore, the objective is only to maximize the total number of relocation requests satisfied. While, in the latter, the aim is to maximize the total profit as the difference between the total revenue associated with the relocation requests satisfied and the total worker cost. In this latter scenario, the worker cost is evaluated multiplying the number of the routes by the unitary personnel cost, i.e., 60 €. The unitary personnel cost has been derived considering a time cost of 12 €/h (as in Boyaci et al., 2013), and that each operator is employed for a duty time of 5 h. While the parameter $rev_r$ (i.e. the revenue associated with the generic request $r$ satisfied) has been set to 20 € like the *penalty cost* for lost customers indicated in the above cited paper.

Instead, the V-AMAT set of instances are characterized by a variable revenue associated with each relocation request. More specifically, $rev_r$, for each relocation request $r$, is made up by two components. The first one represents a variable revenue component (hereafter indicated with VRC), proportional to the actual rent-time of the user that in the simulator generates the request $r$, by a factor equal to 0.29 €/minute. While, the second one is a fixed revenue component (indicated by FRC) associated with the customer satisfaction, i.e., it takes into account that a satisfied customer may require the carsharing service also in the future. This second component, set to 15 €, avoids that the relocation requests associated with customers requiring the carsharing service only for a short rent-time are never satisfied. The value of 15 € arises from the fact that we want to obtain profits comparable to those of AMAT set being the VRC between 1.45 € and 4.35 € since the rent-times of the users are between 5 and 15 minutes.

The next two subsections will show the numerical results. In particular, with reference to the AMAT set, the comparisons between the solutions found by the proposed heuristics and the two MILP formulations will be shown in Subsection 5.1. Moreover, MILP1 and MILP2 will be also compared. While, concerning to the V-AMAT set, the comparisons between the solutions found by the proposed heuristics and the MILP2 are given in Subsection 5.2.

Such comparisons will be discussed not only with reference to the quality of the solution but also regarding the computational times.

### 5.1 Numerical results on the AMAT set

In this section, the experiments carried out on the benchmark instances are described and discussed. Table 3 compares the numerical results obtained by the four heuristic approaches described in Section 3 and Section 4 with the ones reached solving MILP1. The header $\Delta\sigma$ represents the percentage gaps on the number of served requests found by MILP1 ($\sigma^{MILP1}$) and the heuristics ($\sigma^{Heur}$), computed as in the following:

$$\Delta\sigma = \frac{\sigma^{MILP1} - \sigma^{Heur}}{\sigma^{MILP1}} \cdot 100 \ .$$

The header $\Delta K$ represents the gap between the number of workers employed by the heuristics ($K^{Heur}$) and the one employed by the MILP ($K^{MILP1}$), evaluated as in the following:

$$\Delta K = K^{MILP1} - K^{Heur}.$$

The header CPU refers to the computational times required by a heuristic for solving each instance. The structured heuristics yield by far better results than the two greedy ones.

The NNH and MUH do not obtain the optimal solution for eleven instances (with an average $\Delta\sigma$ gap equal to 12.33%) and for ten instances (with an average $\Delta\sigma$ gap equal to 14.97%), respectively; the CH, only for two instances (with an average $\Delta\sigma$ gap equal to 9.54%) and remarkably, the RH always detects the optimum.

As the number of the operators employed is concerned, the NNH uses for nine instances one more worker than the optimal solution while in three instances, one less; in the all other instances, it uses the same number of workers (hereafter, in the paper, it is omitted to indicate the number of instances for which $\Delta K=0$). Moreover, the only three cases with a positive $\Delta K$ (i.e., AMAT 4, AMAT 16 and AMAT 28) are due to the more requests handled by MILP1 than NNH, requiring one more worker.

The MUH, in nineteen instances, uses one more worker than the optimal solution while in two instances, two more. The CH, in ten instances uses one more worker while in four instances, two more and, remarkably, the RH always uses the same number of workers employed in the optimum. RH is suitable to outperform the other heuristics since it generates several feasible solutions and returns the best among them. While, the other heuristics only generate one feasible solution.

Concerning the total computational times, the NNH gives the final solution in 0.04 s, on average, the MUH in 0.06 s, the CH in 0.06 s, and finally, the RH in 3.11 s, against 24.34 s required by the MILP1.

In Table 4, the results obtained considering the E-VReP with the new objective function, introduced in Section 2, and the ones of MILP1 are compared. The header $\Delta\sigma$ denotes the percentage gap between the number of served requests by the E-VReP with the new objective function and the original one. The header $\Delta\pi$ represents the percentage gap between the total profit obtained by the new objective function and the original one, as specified in the following:

$$\Delta\pi = \frac{\pi^{MILP2} - \pi^{MILP1}}{\pi^{MILP1}} \cdot 100.$$

Finally, the header CPU1 and CPU2 show the computational times required by MILP1 and MILP2, respectively.

On eight instances, MILP2 maximizes the total profit compared to MILP1 of 6.5% on average although it decreases the number of served requests of 10.1% on average. Since in MILP2, it could be not convenient (from the profit point of view) to serve the same number of requests handled by MILP1. Moreover, a positive value of $\Delta\pi$ and a negative value of $\Delta K$ always correspond to every negative value of $\Delta\sigma$.

Moreover, it is worth observing that the change applied to the objective function in MILP2 decreases the computational performances since for twenty-three instances the CPU time limit is reached requiring 16,265.4 $s$ on average, against 24.30 $s$ used by MILP1.

In particular, in eight cases, MILP1 handles more requests than MILP2 but using one worker more. Therefore, MILP2 guarantees a higher profit. There is only one case (AMAT 17) in which MILP1 is able to handle more requests than MILP2, using one less worker and then, guaranteeing a higher profit. However, this is a case in which MILP2 reaches the CPU time limit.

Table 5 compares the total profit obtained by the heuristics to the one of MILP2. In particular:

$$\Delta\pi = \frac{\pi^{MILP2} - \pi^{Heur}}{\pi^{MILP2}} \cdot 100$$

and

$$\Delta K = K^{MILP2} - K^{Heur}.$$

The NNH on average decreases the total profit of 18.97% compared to MILP2 while it uses in eleven instances one more worker and, in two instances, two more. The MUH on average decreases the total profit of 32.14% compared to MILP2 while it employs in sixteen instances one more worker, in five instances, two more and finally, in three instances, three more. The CH on average decreases the total profit of 24.77% with respect to MILP2 while it uses, in eleven instances, one more worker, in six instances, two more and in one instance, three more. Finally, the RH on average decreases the total profit only of 1.76% while it uses, only in two instances, one more worker than MILP2.

Concerning the computational times, the NNH on average requires 0.04 $s$, the MUH 0.06 $s$, the CH 0.06 $s$ and finally, the RH 3.16 $s$, against 16,265.4 $s$ of MILP2.

### 5.2 Numerical results on the V-AMAT set

In this subsection, the comparisons between the solutions found by the four proposed heuristics and the ones found by the MILP, on the V-AMAT set, are discussed.

Table 6 compares the total profit obtained by the heuristics to the one of MILP2. The meaning of each column is the same of those of Table 5. The NNH on average decreases the total profit of 16.77% compared to MILP2 while it uses in eight instances one more worker. The MUH on average decreases the total profit of 24.74% compared to MILP2 while it employs in six instances one less worker and in five instances, one more worker. The CH on average decreases the total profit of 30.99% with respect to MILP2 while it uses, in five instances, one fewer worker; in five instances, one more worker; in one instance, two more workers and in one instance, three more workers. Finally, the RH on average decreases the total profit only of 0.61% and only in one instance, it uses one more worker of the MILP2.

Concerning the computational times, the NNH on average requires 0.01 s as well as the MUH while the CH, 0.02 s and finally, the RH 8.36 s, against 8,267.66 s of MILP2.

|  | NNH | | | MUH | | | CH | | | RH | | |
|---|---|---|---|---|---|---|---|---|---|---|---|---|
| Instance | Δσ % | ΔK | CPU | Δσ % | ΔK | CPU | Δσ % | ΔK | CPU | Δσ % | ΔK | CPU |
| AMAT 1 | 0.00 | 0 | 0.04 | 0.00 | 0 | 0.03 | 0.00 | 0 | 0.04 | 0.00 | 0 | 3.94 |
| AMAT 2 | 0.00 | -1 | 0.03 | 0.00 | -1 | 0.03 | 0.00 | -1 | 0.04 | 0.00 | 0 | 2.1 |
| AMAT 3 | 0.00 | -1 | 0.03 | 0.00 | 0 | 0.03 | 0.00 | 0 | 0.03 | 0.00 | 0 | 2.13 |
| AMAT 4 | 14.28 | 1 | 0.04 | 7.14 | 0 | 0.05 | 0.00 | -1 | 0.07 | 0.00 | 0 | 3.4 |
| AMAT 5 | 0.00 | -1 | 0.04 | 0.00 | -1 | 0.88 | 0.00 | 0 | 0.05 | 0.00 | 0 | 3.57 |
| AMAT 6 | 0.00 | -1 | 0.03 | 0.00 | 0 | 0.03 | 0.00 | 0 | 0.03 | 0.00 | 0 | 2.51 |
| AMAT 7 | 0.00 | 0 | 0.03 | 14.28 | -2 | 0.03 | 0.00 | 0 | 0.03 | 0.00 | 0 | 2.49 |
| AMAT 8 | 0.00 | 0 | 0.03 | 0.00 | -1 | 0.04 | 0.00 | -1 | 0.05 | 0.00 | 0 | 2.85 |
| AMAT 9 | 0.00 | 0 | 0.03 | 0.00 | 0 | 0.03 | 0.00 | 0 | 0.03 | 0.00 | 0 | 2.84 |
| AMAT 10 | 9.09 | 0 | 0.27 | 9.09 | -1 | 0.02 | 0.00 | -2 | 0.04 | 0.00 | 0 | 3.15 |
| AMAT 11 | 12.50 | 0 | 0.03 | 0.00 | -1 | 0.03 | 0.00 | 0 | 0.03 | 0.00 | 0 | 2.54 |
| AMAT 12 | 0.00 | 0 | 0.03 | 33.33 | -1 | 0.03 | 0.00 | 0 | 0.03 | 0.00 | 0 | 2.28 |
| AMAT 13 | 0.00 | 0 | 0.04 | 27.27 | -1 | 0.03 | 0.00 | -1 | 0.05 | 0.00 | 0 | 3.09 |
| AMAT 14 | 0.00 | 0 | 0.03 | 0.00 | 0 | 0.03 | 0.00 | 0 | 0.51 | 0.00 | 0 | 2.07 |
| AMAT 15 | 0.00 | 0 | 0.02 | 0.00 | -1 | 0.03 | 0.00 | -1 | 0.03 | 0.00 | 0 | 2.41 |
| AMAT 16 | 18.18 | 1 | 0.03 | 0.00 | 0 | 0.03 | 0.00 | 0 | 0.04 | 0.00 | 0 | 3.35 |
| AMAT 17 | 0.00 | -1 | 0.04 | 0.00 | -1 | 0.04 | 0.00 | -2 | 0.05 | 0.00 | 0 | 5.71 |
| AMAT 18 | 10.00 | 0 | 0.03 | 10.00 | -2 | 0.06 | 0.00 | -1 | 0.07 | 0.00 | 0 | 3.58 |
| AMAT 19 | 0.00 | -1 | 0.03 | 0.00 | 0 | 0.03 | 0.00 | 0 | 0.03 | 0.00 | 0 | 2.74 |
| AMAT 20 | 0.00 | 0 | 0.03 | 0.00 | -1 | 0.03 | 0.00 | -1 | 0.04 | 0.00 | 0 | 3.78 |
| AMAT 21 | 0.00 | -1 | 0.04 | 0.00 | -1 | 0.04 | 0.00 | -1 | 0.05 | 0.00 | 0 | 5.75 |
| AMAT 22 | 14.29 | 0 | 0.03 | 0.00 | -1 | 0.03 | 0.00 | 0 | 0.03 | 0.00 | 0 | 2.07 |
| AMAT 23 | 11.11 | 0 | 0.03 | 11.11 | -1 | 0.03 | 0.00 | -2 | 0.04 | 0.00 | 0 | 2.74 |
| AMAT 24 | 0.00 | 0 | 0.02 | 0.00 | -1 | 0.02 | 0.00 | 0 | 0.03 | 0.00 | 0 | 2.08 |
| AMAT 25 | 0.00 | -1 | 0.03 | 0.00 | -1 | 0.02 | 0.00 | -1 | 0.04 | 0.00 | 0 | 2.06 |
| AMAT 26 | 9.09 | 0 | 0.03 | 9.09 | -1 | 0.03 | 9.09 | 0 | 0.04 | 0.00 | 0 | 3.29 |
| AMAT 27 | 8.34 | -1 | 0.03 | 8.34 | -1 | 0.03 | 0.00 | -1 | 0.04 | 0.00 | 0 | 4.46 |
| AMAT 28 | 20.00 | 1 | 0.03 | 20.00 | 0 | 0.04 | 10.00 | 0 | 0.05 | 0.00 | 0 | 3.86 |
| AMAT 29 | 20.00 | 0 | 0.03 | 0.00 | -1 | 0.03 | 0.00 | -2 | 0.05 | 0.00 | 0 | 3.32 |
| AMAT 30 | 0.00 | 0 | 0.03 | 0.00 | -1 | 0.04 | 0.00 | 0 | 0.05 | 0.00 | 0 | 3.24 |
| **AVERAGE** | **4.90** | | **0.04** | **4.99** | | **0.06** | **0.64** | | **0.06** | **0.00** | | **3.11** |

**Table 2** Comparisons of the four heuristics with MILP1

| Instance | Δσ % | Δπ % | ΔK | CPU1 | CPU2 |
|---|---|---|---|---|---|
| AMAT 1 | 0.0 | 0.0 | 0 | 1.49 | † |
| AMAT 2 | 0.0 | 0.0 | 0 | 0.03 | 779.02 |
| AMAT 3 | 0.0 | 0.0 | 0 | 0.05 | 1,864.60 |
| AMAT 4 | -7.1 | 5.3 | -1 | 123.33 | 3,640.22 |
| AMAT 5 | 0.0 | 0.0 | 0 | 0.13 | 4,200.40 |
| AMAT 6 | -14.3 | 12.5 | -1 | 0.03 | 3,852.14 |
| AMAT 7 | 0.0 | 0.0 | 0 | 0.05 | 1,025.10 |
| AMAT 8 | 0.0 | 0.0 | 0 | 0.68 | 3,765.38 |
| AMAT 9 | 0.0 | 0.0 | 0 | 0.01 | 3,678.24 |
| AMAT 10 | 0.0 | 0.0 | 0 | 2.19 | † |
| AMAT 11 | 0.0 | 0.0 | 0 | 0.06 | 17,049.23 |
| AMAT 12 | 0.0 | 0.0 | 0 | 0.02 | 76.26 |
| AMAT 13 | -9.1 | 6.3 | -1 | 0.41 | 3,658.22 |
| AMAT 14 | 0.0 | 0.0 | 0 | 0.02 | 13,780.46 |
| AMAT 15 | 0.0 | 0.0 | 0 | 0.03 | 13.02 |
| AMAT 16 | -9.1 | 6.3 | -1 | 2.21 | † |
| AMAT 17 | -6.1 | -7.6 | 0 | 562.83 | † |
| AMAT 18 | 0.0 | 0.0 | 0 | 12.75 | † |
| AMAT 19 | 0.0 | 0.0 | 0 | 0.06 | 6,249.01 |
| AMAT 20 | 0.0 | 0.0 | 0 | 0.17 | † |
| AMAT 21 | 0.0 | 0.0 | 0 | 0.4 | † |
| AMAT 22 | 0.0 | 0.0 | 0 | 0.03 | 527.25 |
| AMAT 23 | -11.1 | 8.3 | -1 | 0.12 | 2,529.19 |
| AMAT 24 | 0.0 | 0.0 | 0 | 0.04 | 16,782.88 |
| AMAT 25 | -14.3 | 12.5 | -1 | 0.04 | 4,443.23 |
| AMAT 26 | -9.1 | 6.3 | -1 | 2.65 | † |
| AMAT 27 | 0.0 | 0.0 | 0 | 15.36 | 3,890.55 |
| AMAT 28 | 0.0 | 0.0 | 0 | 3.89 | † |
| AMAT 29 | 0.0 | 0.0 | 0 | 1.05 | 3,658.79 |
| AMAT 30 | -11.1 | 8.3 | -1 | 0.06 | 3,700.04 |
| **AVERAGE** | **-3.04** | **2.27**[2] | | **24.34** | **16,265.40** |

**Table 3** Comparisons between the results of MILP2 and MILP1

---

[2] The average has been computed excluding the instance AMAT 17 for which a negative value of Δπ has been obtained due to the CPU time limit achievement.

|  | NNH | | | MUH | | | CH | | | RH | | |
|---|---|---|---|---|---|---|---|---|---|---|---|---|
| Instance | Δπ % | ΔK | CPU | Δπ % | ΔK | CPU | Δπ % | ΔK | CPU | Δπ % | ΔK | CPU |
| AMAT 1 | 0.00 | 0 | 0.04 | 0.00 | 0 | 0.03 | 0.00 | 0 | 0.04 | 0.00 | 0 | 5.22 |
| AMAT 2 | 27.27 | -1 | 0.03 | 27.27 | -1 | 0.03 | 27.27 | -1 | 0.04 | 0.00 | 0 | 2.1 |
| AMAT 3 | 33.33 | -1 | 0.03 | 0.00 | 0 | 0.03 | 0.00 | 0 | 0.03 | 0.00 | 0 | 2.13 |
| AMAT 4 | 10.00 | 0 | 0.04 | 15.00 | -1 | 0.05 | 20.00 | -2 | 0.07 | 0.00 | 0 | 3.4 |
| AMAT 5 | 18.75 | -1 | 0.04 | 18.75 | -1 | 0.88 | 0.00 | 0 | 0.05 | 0.00 | 0 | 3.57 |
| AMAT 6 | 44.44 | -2 | 0.03 | 11.11 | -1 | 0.03 | 11.11 | -1 | 0.03 | 0.00 | 0 | 2.51 |
| AMAT 7 | 0.00 | 0 | 0.03 | 81.82 | -3 | 0.03 | 0.00 | 0 | 0.03 | 0.00 | 0 | 2.49 |
| AMAT 8 | 0.01 | 0 | 0.03 | 25.01 | -1 | 0.04 | 25.01 | -1 | 0.05 | 0.01 | 0 | 2.85 |
| AMAT 9 | 0.00 | 0 | 0.03 | 0.00 | 0 | 0.03 | 0.00 | 0 | 0.03 | 0.00 | 0 | 2.84 |
| AMAT 10 | 12.50 | 0 | 0.27 | 31.25 | -1 | 0.02 | 37.50 | -2 | 0.04 | 0.00 | 0 | 3.15 |
| AMAT 11 | 15.38 | 0 | 0.03 | 23.08 | -1 | 0.03 | 0.00 | 0 | 0.03 | 0.00 | 0 | 2.54 |
| AMAT 12 | 0.00 | 0 | 0.03 | 77.78 | -1 | 0.03 | 0.00 | 0 | 0.03 | 0.00 | 0 | 2.28 |
| AMAT 13 | 5.88 | -1 | 0.04 | 41.18 | -3 | 0.03 | 23.53 | -2 | 0.05 | 5.88 | -1 | 3.09 |
| AMAT 14 | 0.00 | 0 | 0.03 | 0.00 | 0 | 0.03 | 0.00 | 0 | 0.51 | 0.00 | 0 | 2.07 |
| AMAT 15 | 0.00 | 0 | 0.02 | 42.86 | -1 | 0.03 | 42.86 | -1 | 0.03 | 0.00 | 0 | 2.41 |
| AMAT 16 | 11.76 | 0 | 0.03 | 5.88 | -1 | 0.03 | 5.88 | -1 | 0.04 | 0.00 | 0 | 3.35 |
| AMAT 17 | 5.33 | -1 | 0.04 | 5.33 | -1 | 0.04 | 18.86 | -2 | 0.05 | -8.20 | 0 | 5.71 |
| AMAT 18 | 14.29 | 0 | 0.03 | 57.14 | -2 | 0.06 | 21.43 | -1 | 0.07 | 0.00 | 0 | 3.58 |
| AMAT 19 | 27.27 | -1 | 0.03 | 0.00 | 0 | 0.03 | 0.00 | 0 | 0.03 | 0.00 | 0 | 2.74 |
| AMAT 20 | 0.00 | 0 | 0.03 | 23.08 | -1 | 0.03 | 23.08 | -1 | 0.04 | 0.00 | 0 | 3.78 |
| AMAT 21 | 16.67 | -1 | 0.04 | 16.67 | -1 | 0.04 | 16.67 | -1 | 0.05 | 0.01 | 0 | 5.75 |
| AMAT 22 | 18.18 | 0 | 0.03 | 27.27 | -1 | 0.03 | 0.00 | 0 | 0.03 | 0.00 | 0 | 2.07 |
| AMAT 23 | 23.08 | -1 | 0.03 | 46.15 | -2 | 0.03 | 53.85 | -3 | 0.04 | 0.00 | 0 | 2.74 |
| AMAT 24 | 0.00 | 0 | 0.02 | 27.27 | -1 | 0.02 | 0.00 | 0 | 0.03 | 0.00 | 0 | 2.08 |
| AMAT 25 | 44.44 | -2 | 0.03 | 44.44 | -2 | 0.02 | 44.44 | -2 | 0.04 | 11.11 | -1 | 2.06 |
| AMAT 26 | 17.65 | -1 | 0.03 | 41.18 | -3 | 0.03 | 17.65 | -1 | 0.04 | 0.00 | 0 | 3.29 |
| AMAT 27 | 27.78 | -1 | 0.03 | 33.33 | -2 | 0.03 | 16.67 | -1 | 0.04 | 0.00 | 0 | 4.46 |
| AMAT 28 | 7.14 | 1 | 0.03 | 28.57 | 0 | 0.04 | 14.29 | 0 | 0.05 | 0.00 | 0 | 3.86 |
| AMAT 29 | 28.57 | 0 | 0.03 | 21.43 | -1 | 0.03 | 42.85 | -2 | 0.05 | 0.00 | 0 | 3.32 |
| AMAT 30 | 7.69 | -1 | 0.03 | 30.77 | -2 | 0.04 | 7.69 | -1 | 0.05 | 0.00 | 0 | 3.24 |
| **AVERAGE** | **13.91** | | **0.04** | **26.79** | | **0.06** | **15.69** | | **0.06** | **0.29** | | **3.16** |

**Table 4** Comparisons of the four heuristics with MILP2 in the AMAT set

|  | NNH | | | MUH | | | CH | | | RH | | |
|---|---|---|---|---|---|---|---|---|---|---|---|---|
| Instance | Δπ % | ΔK | CPU | Δπ % | ΔK | CPU | Δπ % | ΔK | CPU | Δπ % | ΔK | CPU |
| V-AMAT 1 | 11.94 | 1 | 0.021 | 11.69 | 1 | 0.026 | 100.00 | 3 | 0.051 | 0.60 | 1 | 15.47 |
| V-AMAT 2 | 7.39 | 0 | 0.004 | 85.03 | 0 | 0.019 | 2.57 | 0 | 0.032 | 0.00 | 0 | 11.7 |
| V-AMAT 3 | 0.30 | 0 | 0.013 | 17.03 | -1 | 0.013 | 17.20 | -1 | 0.023 | 0.05 | 0 | 12.78 |
| V-AMAT 4 | 0.91 | 0 | 0.008 | 17.20 | 0 | 0.008 | 0.98 | 0 | 0.008 | 0.06 | 0 | 6.91 |
| V-AMAT 5 | 100.00 | 1 | 0.002 | 100.00 | 1 | 0.001 | 100.00 | 1 | 0.001 | 4.86 | 0 | 6.17 |
| V-AMAT 6 | 35.66 | 1 | 0.017 | 35.11 | 1 | 0.009 | 24.51 | 1 | 0.019 | 0.92 | 0 | 13.03 |
| V-AMAT 7 | 1.60 | 0 | 0.005 | 1.59 | 0 | 0.004 | 2.00 | 0 | 0.005 | 0.04 | 0 | 4.91 |
| V-AMAT 8 | 16.60 | 0 | 0.015 | 0.79 | 0 | 0.006 | 48.65 | 0 | 0.028 | 0.37 | 0 | 7.09 |
| V-AMAT 9 | 19.30 | 0 | 0.018 | 0.32 | 0 | 0.006 | 20.08 | 0 | 0.014 | 0.00 | 0 | 6.29 |
| V-AMAT 10 | 23.49 | 0 | 0.004 | 24.17 | 0 | 0.004 | 23.67 | 0 | 0.004 | 0.75 | 0 | 7.45 |
| V-AMAT 11 | 18.05 | 1 | 0.017 | 37.42 | -1 | 0.008 | 1.68 | 0 | 0.015 | 0.08 | 0 | 8.02 |
| V-AMAT 12 | 1.50 | 0 | 0.004 | 2.70 | 0 | 0.004 | 1.98 | 0 | 0.003 | 0.14 | 0 | 4.98 |
| V-AMAT 13 | 27.83 | 1 | 0.011 | 12.57 | 0 | 0.012 | 39.16 | 1 | 0.039 | 0.31 | 0 | 12.96 |
| V-AMAT 14 | 1.67 | 0 | 0.005 | 47.50 | 0 | 0.018 | 41.90 | -1 | 0.015 | 0.34 | 0 | 10.31 |
| V-AMAT 15 | 20.10 | 1 | 0.019 | 16.00 | -1 | 0.02 | 24.26 | -1 | 0.055 | 1.01 | 0 | 15.62 |
| V-AMAT 16 | 0.78 | 0 | 0.008 | 28.59 | -1 | 0.007 | 16.44 | 0 | 0.016 | 0.34 | 0 | 10.96 |
| V-AMAT 17 | 25.00 | 1 | 0.016 | 24.64 | 1 | 0.015 | 32.65 | 0 | 0.11 | 1.27 | 0 | 12.68 |
| V-AMAT 18 | 9.90 | 0 | 0.002 | 20.25 | 0 | 0.002 | 20.01 | 0 | 0.001 | 0.02 | 0 | 3.67 |
| V-AMAT 19 | 0.62 | 0 | 0.005 | 23.87 | 0 | 0.018 | 2.01 | 0 | 0.011 | 0.00 | 0 | 9.38 |
| V-AMAT 20 | 1.99 | 0 | 0.004 | 24.43 | 0 | 0.006 | 2.56 | 0 | 0.005 | 0.00 | 0 | 5.13 |
| V-AMAT 21 | 100.00 | 1 | 0.002 | 100.00 | 1 | 0.001 | 100.00 | 1 | 0.001 | 4.66 | 0 | 4.92 |
| V-AMAT 22 | 1.94 | 0 | 0.004 | 1.20 | 0 | 0.004 | 1.94 | 0 | 0.004 | 0.00 | 0 | 4.51 |
| V-AMAT 23 | 14.38 | 0 | 0.005 | 12.90 | -1 | 0.005 | 40.22 | 2 | 0.015 | 0.08 | 0 | 9.29 |
| V-AMAT 24 | 2.75 | 0 | 0.004 | 3.29 | 0 | 0.004 | 2.85 | 0 | 0.003 | 0.00 | 0 | 4.92 |
| V-AMAT 25 | 27.80 | 0 | 0.02 | 38.68 | -1 | 0.01 | 24.82 | -1 | 0.01 | 0.40 | 0 | 8.21 |
| V-AMAT 26 | 3.26 | 0 | 0.003 | 2.55 | 0 | 0.003 | 40.14 | 0 | 0.007 | 1.92 | 0 | 5.9 |
| V-AMAT 27 | 1.56 | 0 | 0.003 | 0.30 | 0 | 0.003 | 32.01 | 0 | 0.006 | 0.00 | 0 | 8.75 |
| V-AMAT 28 | 20.36 | 0 | 0.002 | 21.89 | 0 | 0.001 | 100.00 | 1 | 0.001 | 0.00 | 0 | 7.63 |
| V-AMAT 29 | 4.26 | 0 | 0.002 | 5.49 | 0 | 0.003 | 62.80 | -1 | 0.003 | 0.00 | 0 | 5.77 |
| V-AMAT 30 | 2.17 | 0 | 0.005 | 25.02 | 0 | 0.005 | 2.48 | 0 | 0.005 | 0.15 | 0 | 5.33 |
| **AVERAGE** | **16.77** | | **0.01** | **24.74** | | **0.01** | **30.99** | | **0.02** | **0.61** | | **8.36** |

**Table 5** Comparisons of the four heuristics with MILP2 in the V-AMAT set.

## 6 Sensitivity analysis

In this section, we perform a sensitivity analysis on the main input parameters of the E-VReP, such as the number of relocation requests (instance size) and the revenue associated with each of them.

### 6.1 Sensitivity analysis on the instance size

The sensitivity analysis on the instance size has been carried out considering the results obtained by the exact model MILP2 on both the AMAT and V-AMAT set. The results have been analyzed from the point of view of the profit, of the percentage number of requests served and of the number of operators used.

Concerning the AMAT set, the behavior of the profit, varying the instance size, is shown in Fig. 4. It is possible to observe that the profit increases almost linearly with the number of requests. While, regarding the V-AMAT set, the profit has not such a linear behavior, as shown in Fig. 5. That is coherent with the fact that in this second set of instances the revenue owns a variable component. However, we can observe that until an instance size around 30 requests, the profit is always less than 250 € while, for bigger instance sizes, the profit is almost always greater than 293 €.

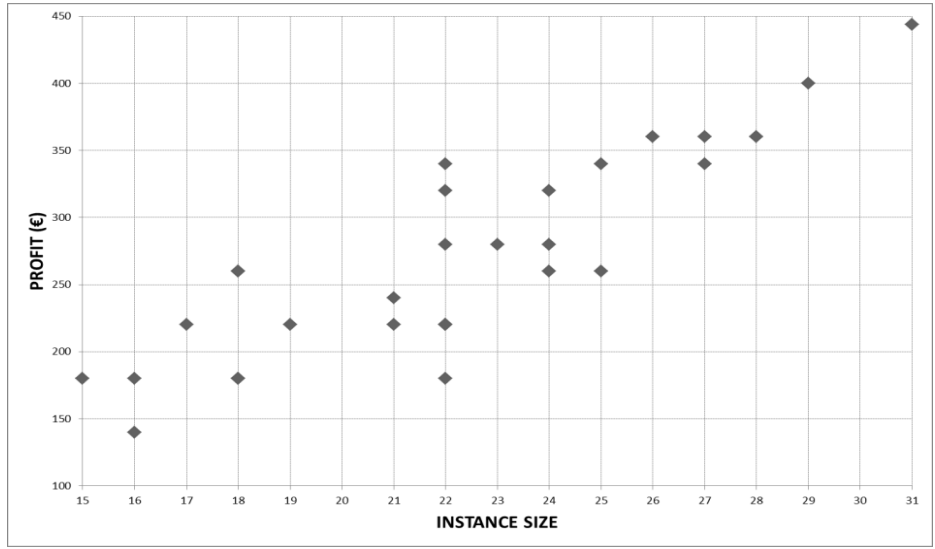

**Fig. 4.** Behavior of the profit obtained by MILP2 varying the instance size in the AMAT set.

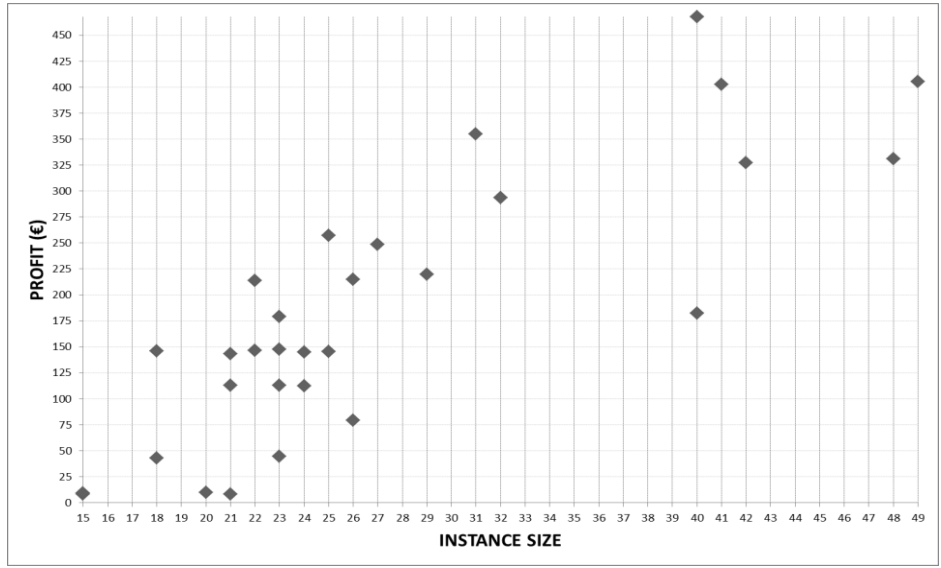

**Fig. 5.** Behavior of the profit obtained by MILP2 varying the instance size in the V-AMAT set.

In Fig. 6 the behavior of the percentage of the requests served varying the instance size for the AMAT set is displayed. It is possible to see that for instances with same size, very different results can be obtained e.g. for the instances with size 22 the percentage of the requests served varies from 55% to 100%. Vice versa the same (or almost the same) percentage of requests served can be obtained for instances with very different sizes (e.g., 80% is obtained both for the smallest instance, with size 15, and for that with size 25). A similar behavior also happens for the V-AMAT set, as shown in Fig. 7. An important difference between the two sets is that while in the first one most instances have a percentage of requests served between 80% and 92%, in the second set the distribution of such percentage is more uniform and in a wider range ([20%, 95%] against [55%,100%] of the former). This behavior may be due to the fact that in the AMAT set the revenues associated with the requests are constant and higher (20 €) than those in the V-AMAT set where can vary between 16.45 € and 19.35 €.

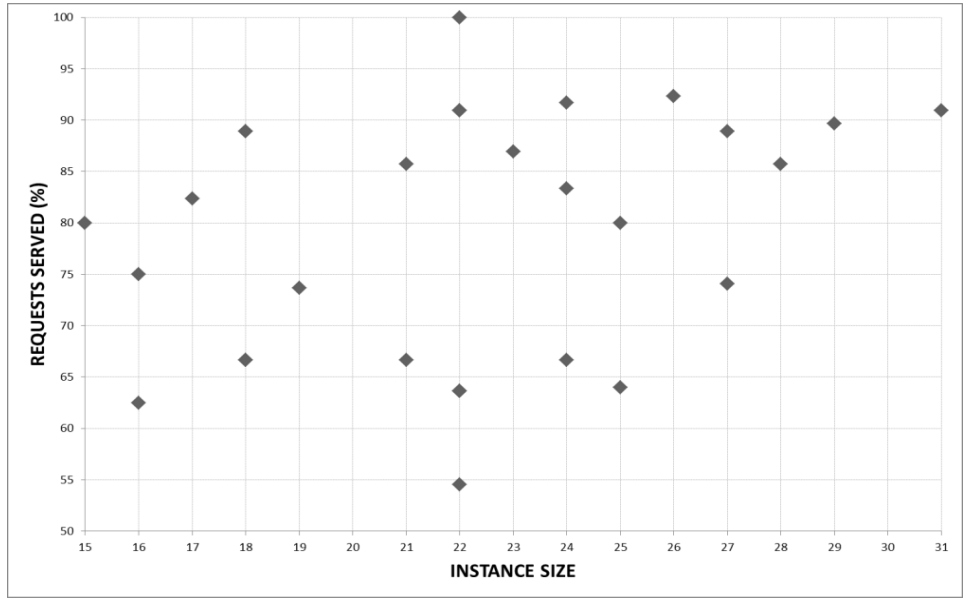

**Fig. 6.** Behavior of the percentage of the requests served varying the instance size in the AMAT set.

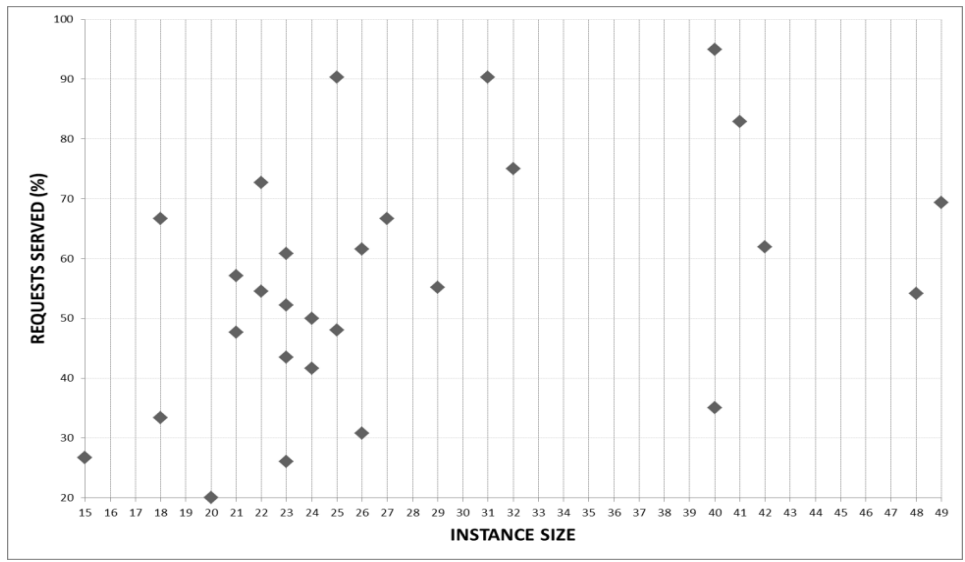

**Fig. 7.** Behavior of the percentage of the requests served varying the instance size in the V-AMAT set.

In Fig. 8, the behavior of the number of the operators used varying the instance size for the AMAT set is displayed. It is possible to observe that for the instances with size between 21 and 27 there is always one instance where the number of operators used is one and one instance where it is two; while, for instances with lower size, only one operator is always used and for those with upper size, two operators.

Concerning the V-AMAT set there is a similar behavior as shown in Fig. 9: for the instance size between 25 and 32 there is always one instance where the number of operators used is one and one instance where it is two; while, for instances with lower size, only one operator is always used and for those with upper size, two or three operators (except for the size 40 for which there is also one instance with one operator).

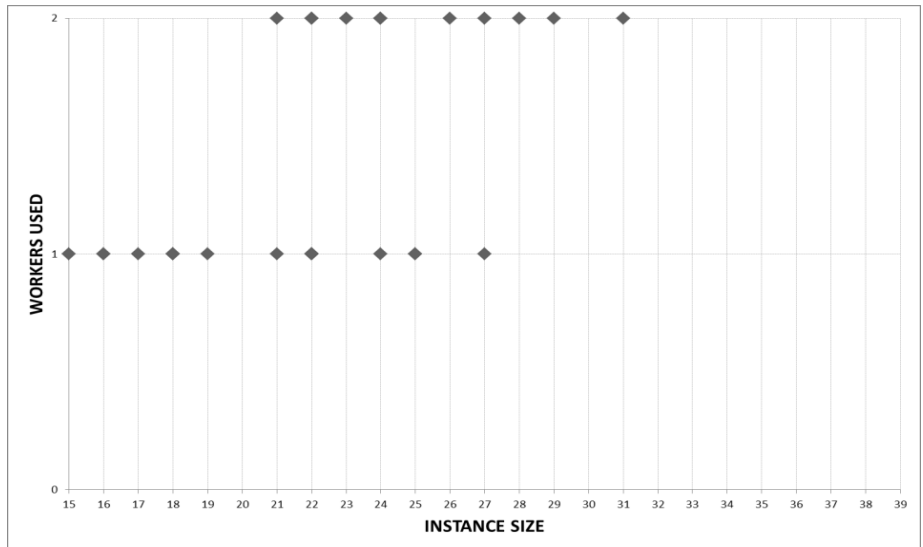

**Fig. 8.** Behavior of the number of workers used varying the instance size in the AMAT set.

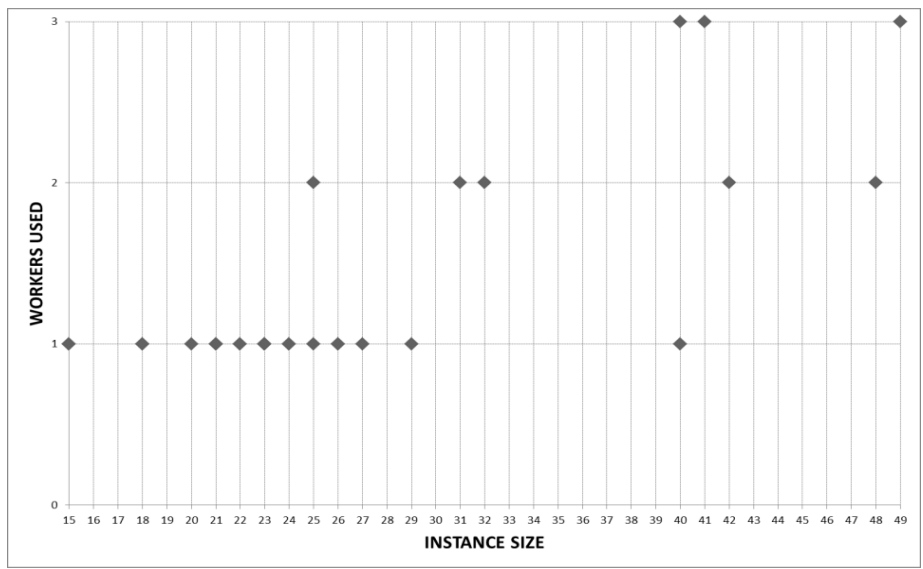

**Fig. 9.** Behavior of the number of workers used varying the instance size in the V-AMAT set.

## 6.2 Sensitivity analysis on the relocation request revenue

With reference to the sensitivity analysis carried out on the revenue associated with each relocation request, we consider the V-AMAT set in order to investigate the impact of the FRC (originally set to 15 € in Subsection 5.2) on the solutions. For this purpose, we tested a significant range of values for FRC making it to vary from 0 to 20 € by 5 €. For each value of FRC, the E-VReP has been solved using the RH heuristic since it gave results comparable with those of the MILP2 (but in by far less CPU time).

Figures 10, 11 and 12 show the behaviors of the profit, of the number of requests served and of the operators used, respectively, for all FRC values. In the plots, a different line is used for each FRC value. Moreover, we do not plot the graphic obtained when FRC = 0, since in this case no relocation is performed being in every instance the total cost (due to the operators) greater than the total revenue.

Fig. 10 shows that the total profit increases almost proportionally as the FRC increases: this behavior is reasonable since the total profit is linear proportional through FRC to the number of requests served and the latter also tends to increase when FRC increases. We observe that when FRC = 5 €, for seven instances the total profit is zero since no relocation is performed, while this happens for four instances when FRC = 10 € and never for greater FRC values. There is only one case where the profit obtained with FRC = 10 dominates the one with FRC = 15 (this is due to the use of a heuristic rather than an exact approach).

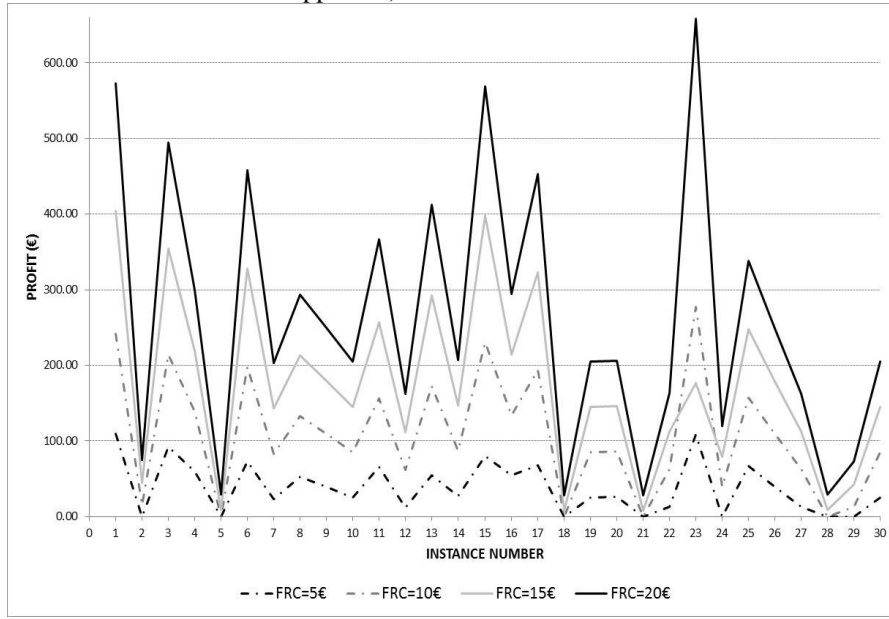

**Fig. 10.** Behavior of the profit varying FRC in the V-AMAT set.

Fig. 11 plots the behaviors concerning the percentage of the number of requests served. As reasonable, increasing FRC, these percentages tend to increase too, but very often they are the same for different FRC values and sometimes (instances 1, 6, 15 and 17) the number of requests served is even greater for FRC = 15 than for FRC = 20. This can be explained by the fact that in the V-AMAT the profit is not merely proportional to the number of the requests served and then it may happen that a fewer number of requests may have a higher total profit (indeed no dominance of the profit is observed).

Finally, Fig. 12 shows the behaviors of the number of workers used and highlights that in general increasing the FRC, more operators need to be employed, but for FRC = 15 and FRC = 20 the number of operators used is the same apart for one instance.

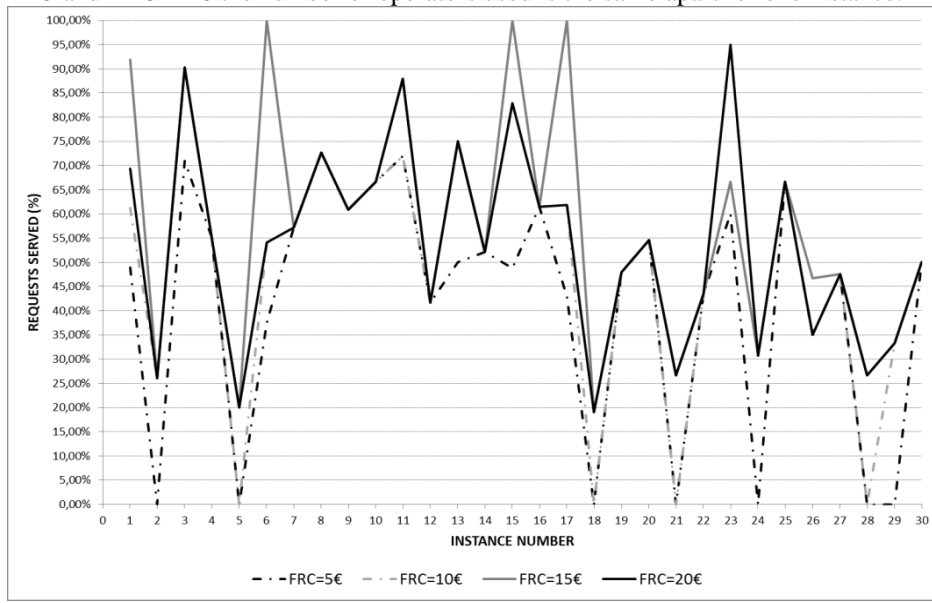

**Fig. 11.** Behavior of the percentage of requests served varying FRC in the V-AMAT set.

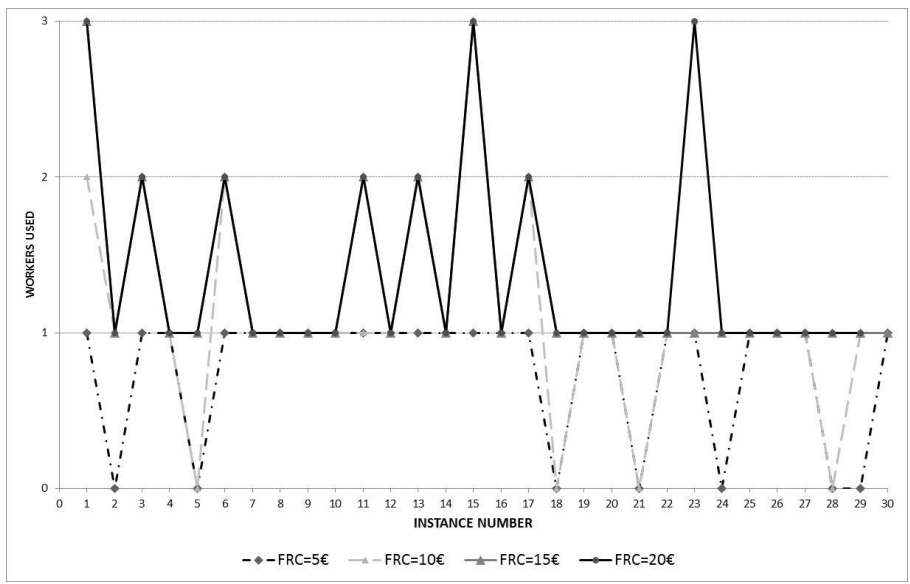

**Fig. 12.** Behavior of the number of the operator used varying FRC in the V-AMAT set.

# 7   Conclusions

The economical sustainability of a new operator based relocation approach for one-way electric carsharing systems (E-VReP), recently proposed in the literature, was investigated in this paper. This aspect was completely neglected till now, despite its importance to understand the practicability of the new relocation approach.

To deal with the economical sustainability of the E-VReP, we assumed that a revenue is associated with each relocation request satisfied and a cost is associated with each worker used. Thus, the original E-VReP objective (given only by the maximization of the relocation requests satisfied) was modified into the maximization of the total profit, i.e., the difference between the total revenue of all the relocation requests satisfied and the total cost of all the workers used.

This new variant of the E-VReP is more challenging than the original problem as empirically proven by the huge rise of the CPU time, required for solving the corresponding MILP (an average rise of about 16,241 *s*, on the AMAT benchmark instances). For this reason, we developed four heuristic approaches, two of which (CH and RH) exploit some general properties of the feasible solutions, specifically proved in this paper.

To test the economical sustainability issue, we built a new set of instances (V-AMAT) through a MATLAB simulator based on realistic data of Milano. In the V-AMAT set, the revenue associated with each relocation request depends on two components: a variable one, proportional to the rent-time associated with the request and a fixed one (FRC). The FRC allows modeling a "future revenue" related to the customer satisfaction since a satisfied customer will ask (likely) for the service also in the future. A sensitivity analysis on the FRC was performed in order to understand its impact on the solutions. Such analysis was possible thanks to the very short CPU time required by the RH heuristic compared to that needed to CPLEX for solving the MILP formulation (on average 8.36 s against 8,267.66 s for each instance, respectively). Despite we used a heuristic approach, we were able to guarantee a high precision as proved by the very small average relative worsening gap between RH and CPLEX (i.e., only 0.61% on the V-AMAT set, when FRC = 15 €).

A significant insight of this work, derived from the above sensitivity analysis, consists in the conclusion that the direct revenue associated with the EV rent by part of the users is not sufficient to cover the worker costs and therefore, it is necessary to include an FRC of at least 15 € to guarantee that relocations are carried out (assuming the values of the data used in our case study).

## References

1. Agenzia Mobilità Ambiente e Territorio, Grafo di offerta stradale versione 0.4 (road network database), (2008) , 2008.
2. A2A. Progetto E-Moving (in Italian), (2013). Available at http://www.e-moving.it/home/cms/emv/.

3.  AMAT. Dati sul traffico veicolare privato sulla rete stradale di Milano – (Origin Destination Matrix, in Italian), (2005). Available at http://www.amat-mi.it/it/downloads/8/.
4.  Archetti, C., Speranza, M., 2004. Vehicle routing in the 1-skip collection problem. Journal of the Operational Research Society 55 , no. 7, 717-727.
5.  Aringhieri, R., Bruglieri, M., Malucelli, F., Nonato, M., 2004. An asymmetric vehicle routing problem arising in the collection and disposal of special waste. Electronic notes in discrete mathematics 17, 41-47.
6.  Ahuja, R. K., Magnanti, T. L., Orlin, J. B., 1993. Network Flows: Theory, Algorithms, and Applications. Prentice Hall.
7.  Barth, M.,Todd, M., 1999. Simulation model performance analysis of a multiple station shared vehicle system. Transportation Research Part C: Emerging Technologies 7, no. 4, 237-259.
8.  Barth, M., Todd, M., Xue, L., 2004. User-based vehicle relocation techniques for multiple-station shared-use vehicle systems. Transportation Research Record 1887, 137-144.
9.  Boyacı, B., Geroliminis, N., Zografos, K., 2013. An optimization framework for the development of efficient one-way car sharing systems. 13th Swiss Transport Research conference.
10. Bruglieri, M., Colorni, A., Luè, A., 2014(a). The vehicle relocation problem for the one-way electric vehicle sharing. Networks, 64(4), 292-305. (ArXive pre-print version of the paper available at: http://arxiv.org/abs/1307.7195v1).
11. Bruglieri, M., Colorni, A., Luè, A., 2014(b). The vehicle relocation problem for the one-way electric vehicle sharing: an application to the Milan case. Procedia Social and Behavioral Sciences 111, 18-27.
12. Correia, G. H. A., Antunes, A. P., 2012. Optimization approach to depot location and trip selection in one-way carsharing systems. Transportation Research Part E: Logistics and Transportation Review 48, 233–247.
13. Correia, G. H. A., Jorge, D. R., Antunes, D. M., 2014. The Added Value of Accounting For Users' Flexibility and Information on the Potential of a Station-Based One-Way Car-Sharing System: An Application in Lisbon, Portugal. Journal of Intelligent Transportation Systems: Technology, Planning, and Operations, 18(3), 299-308.
14. Cucu T., Ion, L., Ducq, Y. and Boussier, J.-M, 2009. Management of a public transportation service: Carsharing service. Proceedings of The 6th International Conference on Theory and Practice in Performance Measurement and Management.
15. Di Febbraro, A., Sacco, N., Saeednia, M., 2012. One-way carsharing: Solving the relocation problem. Transportation Research Board 91st Annual Meeting.
16. Du, Y., Hall, R., 1997. Fleet sizing and empty equipment redistribution for center-terminal transportation networks. Management Science 43, no. 2, 145-157.
17. European Commission, *White paper: Roadmap to a single european transport area–towards a competitive and resource efficient transport system*, (2011). *COM (2011)*, vol. 144.

18. Fourer, R., Gay, D., Kernighan, B. W., 2002. The ampl book. Duxbury Press, Pacific Grove.
19. Hafez, N., Parent, M., Proth, J. M., 2011. Managing a pool of self service cars. Intelligent Transportation Systems. Proceedings. 2001 IEEE, IEEE, 943-948.
20. Jorge, D., Correia, G. H. A., 2013. Carsharing systems demand estimation and defined operations: a literature review. European Journal of Transport and Infrastructure Research, 13, 201-220.
21. Jung, J., Chow, J. Y., Jayakrishnan, R., & Park, J. Y., 2014. Stochastic dynamic itinerary interception refueling location problem with queue delay for electric taxi charging stations. Transportation Research Part C: Emerging Technologies, 40, 123-142.
22. Kek, A. G., Cheu, R. L., Meng, Q., Fung, C. H., 2009. A decision support system for vehicle relocation operations in carsharing systems. Transportation Research Part E: Logistics and Transportation Review 45, no. 1, 149-158.
23. Nourinejad, M., Roorda, M. J., 2014. A dynamic carsharing decision support system. Transportation Research Part E, 66, 36--50.
24. Touati-Moungla, N., Jost, V., 2012. Combinatorial optimization for electric vehicles management. Journal of Energy and Power Engineering 6, no. 5, 738-743.
25. Wang, H., Cheu, R. , Lee, D. H., 2010. Logistical Inventory Approach in Forecasting and Relocating Share-use Vehicles. Advanced Computer Control (ICACC), 2010 2nd International Conference on , 5, 314-318.